\documentclass[11pt]{article}
 \usepackage{amsfonts}
 \usepackage{amssymb}

\newtheorem{theorem}{Theorem}

\newcommand\bqa {\begin{eqnarray}}
\newcommand\eqa {\end{eqnarray}}
\newcommand{\beq}{\begin{eqnarray}}
\newcommand{\eeq}{\end{eqnarray}}
\newcommand{\be}{\begin{array}}
\newcommand{\ee}{\end{array}}

\newcommand{\pr}{\partial}
\newcommand\apr {\overline {\partial }}

\newcommand{\E}{{\cal E}}
\newcommand{\s}{{\cal S}}

\newcommand{\D}{{\cal D}}

\newcommand{\RR}{{\mathbb R}}

\newcommand{\adot}{\dot{\alpha}}
\newcommand{\bdot}{\dot{\beta}}

 \newcommand{\g}{{\mathfrak g}}

\begin{document}

\def\t{\theta}
\def\T{\Theta}
\def\w{\omega}
\def\ov{\overline}
\def\a{\alpha}
\def\b{\beta}
\def\g{\gamma}
\def\s{\sigma}
\def\l{\lambda}
\def\wt{\widetilde}

\def\P{\mathbb{P}}
\def\O{\mathcal{O}}
\def\J{\mathcal{J}}
\def\N{\mathcal{N}}
\def\H{\mathcal{H}}

\begin{flushright}
 ITEP-TH 84-2000 \\
\end{flushright}

\vskip 10mm

\begin{center}
{ \bf \large Harmonic Twistor Formalism  and Transgression
              on Hyperk\"{a}hler manifolds.}
\end{center}

\smallskip

\begin{center}
A.Gerasimov, A.Kotov
\end{center}

\centerline{\sl ITEP, B. Cheriomushkinskaya 25, Moscow, 117259 Russia}

\vskip 10mm

\smallskip





In this paper we continue our study of the transgression of
characteristic classes of hyperholomorphic bundles on
hyperk\"ahler manifolds \cite{GK}. In the previous paper the
global construction for the fourth order transgression of the
Chern character form on a compact hyperk\"ahler manifold was
proposed. In addition, the explicit expression for the
transgression of the Chern character arising in the application of
the local families index theorem was found. This construction was
local over the base of the fibration. It is natural to look for a
 local derivation  of the transgression of the Chern character
forms for an arbitrary hyperholomorphic bundle. In this paper we
give the general local construction for an arbitrary
hyperholomorphic bundle on a four-dimensional dimensional
hyperk\"ahler manifold. Note that in $d=4$ the condition on a
hermitian bundle to be hyperholomorphic is equivalent to the
(anti)self-duality condition on the corresponding connection. We
propose an explicit expression for the fourth order transgression
$T(\E)$  of the top degree part of the Chern character form for an
 arbitrary vector bundle $\E$ supplied with a self-dual connection.
The construction is local and thus is applicable to an arbitrary
four dimensional hyperk\"ahler manifold $M$. Locally the Chern
character form is exact and we have:

$$ch_{[2]}(\mathcal{E})=vol_M \Delta^2T $$ for a volume form
$vol_M$.  Remarkably, the explicit expression for $T(\E)$ is
non-trivial even for a linear bundle $\E$.

In our derivation we essentially use the harmonic twistor
approach, the variant of the twistor formalism developed in
\cite{R,GIKOS1,GIKOS2}. In  twistor approach \cite{WW}  one  codes
the information about self-dual connections on a vector bundle in
terms of  holomorphic structures on a bundle over the twistor
fibration $Z_M\rightarrow M$ with a fiber being $S^2$. Remarkably,
 the proposed expression for $T(\E)$ is given in terms of the
determinant of the $\overline{\partial}_A$-operator in the sense
of Quillen \cite{Q} acting on sections of the holomorphic  bundle
restricted to the fibers. This implies that the results of this
paper may be connected with the local families index for the
twistor fibration.

Rather straightforwardly, the construction described in this paper
 may be generalized to  hyperk\"ahler manifolds of an arbitrary dimension.
We are going to discuss the general construction connecting the
approaches of this paper and of \cite{GK} in the future
publication.

\vskip 15mm

{\bf Acknowledgements}: During the course of this work, the
authors benefited from helpful conversations with A.Levin,
V.Rubtsov. The work of A.G. was partially supported by RFBR grant
98-01-00328 and Grant for the Support of Scientific Schools
00-15-96557. The work of A.K. was partially supported by RFBR
grant 99-01-01254 and Grant for the Support of Scientific Schools
00-15-99296.

\vskip 4mm
\begin{center}
\large \bf Harmonic twistor formalism.
\end{center}
\vskip 4mm

In this section we give a short account of the harmonic formalism
closely following  the presentation given in \cite{GIKOS1,GIKOS2}.
Let $M$ be a $4-$dimensional manifold. Holonomy group in four
dimensions is a tensor product $Sp(1)\otimes Sp(1)$ and $T^* M$
naturally splits: $$ T^* M=\H_L\otimes \H_R, $$ where $\H_{L,R}$
are
 $Sp(1)$-bundles over $M$ with the connection forms
$\omega^{\dot{\alpha}}_{(L,R) \dot{\beta}}$  such that $$ \nabla
h^{\dot{\alpha}}=h^{\dot{\beta}}\omega^{\dot{\alpha}}_{L\dot{\beta}},
\hspace{4mm} \nabla e^{\alpha}=e^{\beta}\omega^{\alpha}_{R\beta}.
$$ In this notations
$\theta^{\dot{\alpha}\alpha}=h^{\dot{\alpha}}\otimes e^{\alpha}$
is the basis of $1-$forms.

Let $Z$ be a total space of $\H_L\setminus 0$. Let us introduce
harmonic variables $u^{\pm\adot}$, $u^{\pm}_{\adot}$ with the
following properties \bqa \label{updown1}
 u^{+}_{\adot}=u^{+ \bdot}\epsilon_{\adot \bdot}, \hspace{4mm}
 u^{+ \adot}=u^{+}_{\bdot}\epsilon^{\adot \bdot} \\ \label{updown2}
u^{-}_{\adot}=u^{- \bdot}\epsilon_{\bdot \adot} , \hspace{4mm}
u^{- \adot}=u^{-}_{\bdot}\epsilon^{\bdot \adot}
 \eqa
where $\epsilon^{\adot \bdot}=-\epsilon_{\adot
\bdot}=\epsilon_{\bdot \adot}$.

They define a frame in the bundle $\H_L$. We will consider the
spherical bundle defined by the condition: \bqa \label{norm}
 u^+_{\adot}u^{+ \adot}=0 , \hspace{4mm}
 u^+_{\adot}u^{-\adot}=1
 \eqa
supplied with the reality condition: \bqa
 \overline{u^{\pm \adot}}=u^{\pm}_{\adot}
\eqa

Let $M$ be a hyperk\"ahler manifold. Then one could chose the
trivial connection on $\H_L$
($\omega^{\dot{\alpha}}_{L\dot{\beta}}=0$). Let us introduce the
bases of vertical forms $\t^{\pm\adot}$ and horizontal forms
$\eta^{\pm\a}$ with respect to projection $Z\rightarrow M$ as
follows: $$ \t^{\pm\adot}=du^{\pm\adot}, \hspace{4mm}
\eta^{\pm\a}=u^{\pm\adot}\theta^{\dot{\alpha}\alpha}. $$

The variables $u^{\pm\adot}$ parameterize  complex structures on
$T^*_x M$ compatible with the hyperk\"ahler structure. At each
point $(x,u)\in Z$  the forms $\eta^{+\a}(\eta^{-\a})$ span the
distribution of the holomorphic (antiholomorphic) forms with
respect to the complex structure $u$ on $T^*_x M$.

We will be interested in the local properties of the self-dual
connections on vector bundles. Thus we could consider the flat
space $M=\RR^4$ with the standard metric $g_{\mu \nu}=\delta_{\mu
\nu}$ in coordinates $x^{\mu}$.  As  the sections of  $Sp(k)\times
Sp(1)-$ bundle, the coordinates $x^{\mu}$ may be written as
$x^{\a\adot}$ with the reality condition
$x^{\a\adot}=\overline{x_{\a\adot}}$, where \bqa
x_{\a\adot}=\epsilon_{\a\b}\epsilon_{\adot\bdot}x^{\b\bdot}=
\epsilon_{\b\a}\epsilon_{\bdot\adot}x^{\b\bdot} \label{updown3}\\
x^{\a\adot}=\epsilon_{\a\b}\epsilon_{\adot\bdot}x_{\b\bdot}=
\epsilon_{\b\a}\epsilon_{\bdot\adot}x_{\b\bdot}. \label{updown4}
\eqa For the complex structure defined by $u$ the
(anti)holomorphic coordinates are:\bqa x^{\pm
\a}=u^{\pm}_{\adot}x^{\a\adot} \\ x^{\pm}_{\a}=u^{\pm
\adot}x_{\a\adot} \eqa Now the reality condition can be written as
follows: $\overline{x^{\pm}_{\a}}=x^{\mp \a}$. In fact $$
\overline{x^{\pm}_{\a}}=\overline{u^{\pm
\adot}}\overline{x_{\a\adot}}= u^{\mp}_{\adot}x^{\a\adot}=x^{\mp
\a} $$ Using (\ref{updown1}, \ref{updown2}, \ref{updown3},
\ref{updown4}), one gets \bqa x^{+ \a}=x^{+}_{\b}\epsilon^{\b\a},
\hspace{4mm} x^{+}_{\a}=x^{+ \b}\epsilon_{\b\a} \\ x^{-
\a}=x^{-}_{\b}\epsilon^{\a\b}, \hspace{4mm} x^{-}_{\a}=x^{-
\b}\epsilon^{\a\b} \eqa

For the differential operators $\partial^{\pm}_{\a}\equiv u^{\pm
\bdot}\partial_{\a \bdot}$ we have the following simple relation:

 \bqa
\partial^{\pm}_{\a}=\frac{\partial}{\partial x^{\mp \a}},
\eqa

Taking into account the normalization conditions (\ref{norm}), we
get the expression for the Laplace operator: \bqa \triangle \equiv
\pr^{\mu}\pr_{\mu}=2
\partial^{\pm}_{\a}\partial^{\mp \a} \label{laplacian}
\eqa

In the following we will use the realization of $sl(2)$ algebra:
\bqa
 &\left[D^{++};D^{--}\right]&=D^0 \\
  &\left[ D^0 ; D^{\pm \pm}\right]&=\pm 2D^{\pm \pm} \eqa

 by the first order differential operators: \bqa
 D^{++}&=&u^{+\adot}\frac{\pr}{\pr u^{-\adot}}\\
 D^{--}&=&u^{-\adot}\frac{\pr}{\pr u^{+\adot}}\\
D^{0}&=&u^{+\adot}\frac{\pr}{\pr
u^{+\adot}}-u^{-\adot}\frac{\pr}{\pr u^{-\adot}}\eqa

with the properties:  \bqa D^{++}u^{+\adot}&=&0 \\
D^{++}u^{-\adot}&=&u^{+\adot} \eqa

One could introduce the formal analog of integration of the
functions of the harmonic variables as follows. Let $f^{(q)}(u)$
be a function  of the charge $q$ ($D^0(f^{(q)})=qf^{(q)}$). Then
it has the following expansion: \bqa
 f^{(q)}(x;u^{\pm})= \sum f^{(\adot_1,\cdots \adot_{n+q}\bdot1
 \cdots  \bdot_n)}
 u^+_{\adot_1}\cdots u^+_{\adot_{n+q}}u^-_{\bdot_1}\cdots
 u^-_{\bdot _n}
 \eqa
 The integration may be defined by the conditions: \bqa \label{integr}
 &\int d^2u \, 1&=1 \\
 &\int d^2u u^+_{(i_1\cdots i_n  }u^-_{ j_1 \cdots j_m
 )}&=0 \, n+m>0 \eqa
Note that the integral of the function with non-zero charge is
zero. Integration rules defined by (\ref{integr}) have the usual
property to be zero for a total derivative: \bqa \int d^2u
D^{++}f=0 \eqa

Thus defined integration is equivalent to the usual integration
over the sphere $S^2$ in terms of coordinates $u^{\pm,\adot}$ but
has a virtue to be defined algebraicaly.

Consider a hermitian vector bundle $\E$ on $M$ with a connection
$\nabla$. For the holomorphic structure defined by $u$ the
holomorphic and antihiolomorphic parts of the connection are given
by: \bqa \label{connect}
\nabla^{\pm}_{\a}=u^{\pm\adot}\nabla_{\a\adot}, \hspace{4mm}
A^{\pm}_{\a}=u^{\pm\adot}A_{\a\adot} \eqa

The curvature of the connection $\nabla$ has the representation:
$$ F_{\a\adot
\b\bdot}=\epsilon_{\adot\bdot}f_{\a\b}+\Omega_{\a\b}f_{\adot\bdot},
$$ where $\epsilon$, $\Omega$ are antisymmetric, $f$ is symmetric.

Let the connection $\nabla$ on  $\E$ be self-dual. Thus the
self-dual part of the curvature $f_{\adot\bdot}$ is zero and we
have: \bqa \label{spincurv}
F_{\a\adot\b\bdot}=\epsilon_{\adot\bdot}f_{\a\b} \eqa

From the conditions (\ref{norm}) it follows that the connection is
integrable on the holomorphic and antiholomorphic hyperplanes:
\bqa [\nabla^{\pm}_{\a},\nabla^{\pm}_{\b}]&=&0,\\ {[}
\nabla^{+}_{\alpha},\nabla^{-}_{\beta} {]}&=&f_{\a\b} \eqa

The operators $\nabla^{\pm}_{\adot}$ which lead to self-dual
connection may be characterized by the following set of equations:
 \bqa \label{conect1}
 \left[\nabla^{\pm}_{\a};\nabla^{\pm}_{\b}\right]&=&0\\ \label{conect2}
 \left[D^{\pm \pm};\nabla^{\pm}_{\a}\right]&=&0 \\ \label{conect3}
 \left[D^0;\nabla^{\pm}_{\a}\right]&=&\pm\nabla^{\pm}_{\a}\eqa

The main advantage of the representation
(\ref{conect1})(\ref{conect2})(\ref{conect3}) is the possibility
to use the $u$-dependent gauge transformations for finding
explicit solutions of the self-duality conditions. This goes as
follows. Locally the first equation (\ref{conect1}) allows to
represent the positive part of the harmonic connection as the pure
gauge with zero-charge gauge parameter $U$: \bqa
\nabla^+_{\a}&=&U^{-1}\partial^+_{\a} U
\\
 q(\nabla^+_{\a})&=&1 \Rightarrow q(U(x,u))=0 \eqa

 After the gauge transformation with parameter
$U(x,u)^{-1}$ we get the set of covariant derivatives: \bqa
\nabla^+_{\a}&=&\partial^+_{\a} \\
 \D^{++}&=&D^{++}+V^{++}=D^{++}+UD^{++}U^{-1} \\
 \D^0&=&D^{0} \eqa

Now the only constraint on the function $V^{++}$ with $q=2$ comes
from (\ref{conect2}):\bqa \label{sol}
 \frac{\partial}{\partial x^{-\a}} V^{++}=0 \eqa
 The solution of this equation obviously is given by an arbitrary
 function of $X^{+\a}$, $u^{\pm \adot}$ with the total charge $q=2$.
Taking into account the properties of the integral (\ref{integr}) we could
reconstruct the gauge field from the solution of the equation
(\ref{sol}):\bqa
 A_{\a \adot}=\int d^2u u^-_{\adot}(U^{-1}\pr^+_{\a}U) \eqa

As an example consider the following  matrix-valued function
corresponding to the gauge group Sp(1) \cite{Inst1}\cite{Inst2}:
\bqa (V^{++})^j_i = x^{+j}x^+_i/\rho^2 \eqa This leads to: \bqa
(U)_i^j&=&(1+\frac{x^2}{\rho^2})^{-\frac{1}{2}}(\delta_i^j+x^{+j}x^-_i/\rho^2)\\
A_{\a \adot i}^{j}&=&\frac{1}{\rho^2+x^2}(\frac{1}{2}x_{\a
\adot}\delta^{j}_i+\epsilon_{i\a}x^j_{\adot}) \eqa Thus we get the
one-instanton solution \cite{BPST} with the center at $x=0$ and
the size $\rho$.

\vskip 4mm
\begin{center}
\large \bf Forth order transgression of the second Chern class
\end{center}
\vskip 4mm

According to the general considerations in \cite{GK} it is natural
to expect that locally the Chern character form of a
hyperholomorphic bundle $\E$  over a hyperk\"ahler manifold $M$
admits  the forth order transgression: \bqa \label{transgr}
ch(\mathcal{E})=dd_Id_Jd_K(\tau(\E)) \eqa where
$d_I=IdI^{-1},d_J=JdJ^{-1},d_K=KdK^{-1}$ are exterior derivative
operators twisted by the compatible complex structures $I,J,K$.
For four-dimensional hyperk\"ahler manifold this relation
simplifies: \bqa \label{transgr4} ch_{[2]}(\mathcal{E})=vol_M
\Delta^2 T(\E) \eqa Here $ch_{[2]}$ is a degree four component of
the Chern character, $\, vol_M$ is the volume form on $M$ and
$\Delta$ is the Laplace operator.

In this paper we prove the relation (\ref{transgr4}) using the
harmonic twistor formalism and give the representation for $T$ in
terms of the determinant of the first order differential operator:
\begin{theorem}\label{th}
{Let $\E$ be a hermitian vector bundle  on a hyperk\"ahler four
dimensional manifold $M$ with a self-dual connection $\nabla$ and
the curvature form $F=\nabla^2$. The following local expression
for the top degree part of the Chern character form holds: \bqa
\label{conject} ch_{[2]}(\E)&=-&\frac{1}{8\pi^2}Tr(F\wedge F)=vol_M
\Delta^2 T(\E)
\\
 T(\E)&=&\frac{1}{16\pi^2}\log(Det(D^{++}+V^{++})/Det(D^{++}))\eqa}
\end{theorem}

 The determinant here is essentially the determinant of the
operator  $\apr$ ( in holomorphic parameterization we have
$D^{++}\sim\pr_{\overline{z}}$).

In our parameterization the formulas (\ref{var5}) (\ref{faltap})
from the Appendix take the form: \bqa \label{variat}
 \delta log Det (D^{++}+V^{++})=\frac{1}{2}\int d^2u (\delta V^{++} V^{--}(V^{++}))
 \eqa
 with the condition:\bqa
 D^{++}V^{--}-D^{--}V^{++}+[V^{++},V^{--}]=0. \label{flat}\eqa

Locally one could represent $A^{+}$-component as a pure gauge:
 \bqa
 &A^+_\a &= U^{-1}\partial^+_\a (U) \
\eqa In this  parameterization we have: \bqa
V^{++} &=& 0 \\
 V^{--} &=& 0 \\
 \nabla^-_{\a}&=& [D^{--}, \nabla^+_{\a}]=U^{-1}[\D^{--}, \partial^+_\a]U \\
 A^-_\a &=& -U^{-1}\partial^+_\a (V^{--})U +U^{-1}\pr^-_{\a}(U)\\
 f_{\a \b}&=&  U^{-1}[\partial^+_\a , \partial^-_\b -\partial^+_\b (V^{--})]U =-
U^{-1}\partial^+_\a\partial^+_\b (V^{--})U
 \eqa
 Making gauge transformation with the gauge parameter $U^{-1}$, we
 get the following  representation:\bqa
A^+_\a&=&0 \\
 V^{++}&=&-D^{++}(U)U^{-1} \\
 V^{--}&=&-D^{--}(U)U^{-1} \\
 A^-_\a&=&- \pr^+_\a (V^{--}) \\
 f_{\a
 \b}&=&-\pr^+_\a\pr^+_\b V^{--} \label{curvature}
 \eqa

The identity (\ref{variat}) together  with (\ref{laplacian}) gives
us: $$ \frac{1}{2}\triangle^2 log Det(D^{++}+V^{++})=
\partial^{+\a}\partial^-_{\a}\partial^{\b\bdot}\partial_{\b\bdot}log
Det (D^{++}+V^{++})=$$ $$ =\frac{1}{2}\int \partial^{+\a}\partial^{-}_{\a}
\pr^{\b\bdot} tr V^{--}\pr_{\b\bdot}(V^{++}) $$

Since $\pr_{\b\bdot}=u^+_{\bdot}\pr^-_{\b}+u^-_{\bdot}\pr^+_{\b}$
and $\pr^{+\b}(V^{++})=0$ we obtain $$\triangle^2 log
Det(D^{++}+V^{++}) =\int \pr^-_{\a}tr \pr^{+\a} \pr^{+\b}(V^{--})
\pr^-_{\b}(V^{++}) $$ Now we have: $$ \pr^-_{\b}(V^{++})=[D^{--},
\pr^+_{\b}]V^{++}= -\pr^+_{\b} D^{--}(V^{++}). $$ Using the
flatness condition (\ref{flat}), one could replace
$D^{--}(V^{++})$ by $\D^{++}(V^{--})$:
$$\D^{++}(V^{--})=D^{++}(V^{--})+[V^{++}, V^{--}].$$ Therefore we
get: $$ \pr^-_{\b}(V^{++})=-\pr^+_{\b} \D^{++}(V^{--})=
-\D^{++}\pr^+_{\b}(V^{--}). $$ Taking into account
(\ref{curvature}), one derives:$$ \triangle^2 log
Det(D^{++}+V^{++})=\int \pr^{-\a}
 tr f^{\a\b} \D^{++}\pr^+_{\b}(V^{--})
$$ Now from the Bianchi identity $\nabla^{-}_{\a}(f^{\a\b})=
\pr_{\a}^- (f^{\a\b}) +[A_{\a}^- ,f^{\a\b}]=0$ one can obtain $$
\triangle^2 log Det(D^{++}+V^{++})=\int tr f^{\a\b}
\nabla^{-}_{\a}\D^{++}\pr^+_{\b}(V^{--})= $$ $$ =\int tr f^{\a\b}
[\nabla^{-}_{\a},\D^{++}]\pr^+_{\b}(V^{--})+ tr
f^{\a\b}\D^{++}\nabla^{-}_{\a}\pr^+_{\b}(V^{--})= $$ $$ =\int tr
f^{\a\b}  f_{\a\b} - tr \D^{++}(f^{\a\b})f_{\a\b} + D^{++}(tr
f^{\a\b}  f_{\a\b}). $$ The second and the third terms are zero
 due to the relations $\D^{++}(f^{\a\b})=0$ and  \\
 $\int D^{++}(tr f^{\a\b} f_{\a\b})=0$.
Therefore we get the simple identity \bqa \triangle^2 log
Det(D^{++}+V^{++})=\int tr f^{\a\b} f_{\a\b} \eqa Taking into
account the relation $(F\wedge F)=-\frac{1}{2}f_{\alpha \beta}f^{\alpha
\beta}vol_M$ we  complete the proof of the theorem.

\vskip 4mm
\begin{center}
\large \bf Explicit calculation for one-instanton connection.
\end{center}
\vskip 4mm

There is a  well known explicit formula for the density of the
topological charge of the gauge field describing one-instanton
solution of the self-duality equations: \bqa 2trF\wedge F=-\Delta^2
\log(1+\frac{x^2}{\rho^2})vol_M \eqa Here the instanton with the
center at $x=0$ has the size $\rho$. This formula is obviously a
particular case of our general formula (\ref{conject}) with: \bqa
\label{t'Hooft}
 16\pi^2T=\log (1+\frac{x^2}{\rho^2}) \eqa

In this section we show how our general expression for $T$ reduces
to (\ref{t'Hooft}). Consider the expansion of the determinant:
\bqa \label{expans}\log Det(1+\frac{1}{D^{++}}V^{++})= \sum_{k=1}
(-1)^k\frac{1}{k}\int d^2u Tr(\frac{1}{D^{++}}V^{++})^k \eqa

Taking into account the simple identity: \bqa
 D^{++}(x^+_ix^{-j})=x^+_ix^{+j} \eqa

let us analyze first terms of the expansion: \bqa \label{iterat}
\frac{1}{D^{++}}&x^+_ix^{+j}/\rho^2&+\frac{1}{D^{++}}x^+_ix^{+j}/\rho^2\frac{1}{D^{++}}x^+_ix^{+j}/\rho^2+\cdots
= \\
=&x^+_ix^{-j}/\rho^2&+\frac{1}{D^{++}}x^+_ix^{+l}x^-_lx^{+j}/\rho^4+\cdots
= \\ = &x^+_ix^{-j}/\rho^2&+x^+_ix^{-j}|x|^2/\rho^4+.. \eqa

Here we have used the relations: \bqa x^+_jx^{+j}&=&0 \\
x^+_jx^{-j}&=&|x|^2 \eqa

It is clear that different terms in the expansion are connected by
simple relations. Taking the integral over $u$ variable  we get
for the full series: \bqa
 16\pi^2T=\sum (-1)^k\frac{1}{k} \frac{|x|^{2k}}{\rho^{2k}}
 =\log(1+\frac{|x|^2}{\rho^2}) \eqa

\vskip 4mm
\begin{center}
\large \bf Appendix: Determinants of Cauchy-Riemann operators over
Riemann surface
\end{center}
\vskip 4mm

In this section we recall basic facts about the determinants of
$\apr_A$-operators and prove the identities used in the main body
of the paper.

Let $M$ be a compact $1-$dimensional complex manifold and $E$ be a
smooth vector bundle over $M$. Let $\nabla_A$ be a holomorphic
connection. We denote $(1,0)$ and  $(0,1)$ components of
$\nabla_A$ as $\pr_A$ and $\bar{\pr_A}$ respectively and identify
the affine space $\mathcal{A}$ of $\bar{\pr_A}-$operators with the
space of holomorphic structures in $E$.

Let $\triangle_A$ be a Laplace operator  written as follows:
$\triangle_A=\bar{\pr_A}^*\bar{\pr_A}$. Assume that $\bar{\pr_A}$
is invertible and $\pr_A$ is conjugated to $\bar{\pr_A}$with
respect to a suitable hermitian metric.

\begin{theorem}(Quillen, \cite{Q})
{Let $A_0$ be a base point. Then there exists a unique up to a constant holomorphic function
$Det(A_0,A)$ on $\mathcal{A}$ such that
\bqa
Det\triangle_A=e^{-||A-A_0||^2}|Det(A_0,A)|^2. \label{quillen}
\eqa
Here  $Det\triangle_A =exp(-\frac{\pr}{\pr s}_{|s=0}Tr\triangle_A^{-s})$
is $\zeta -$regularalized determinant of $\triangle_A$ and
$$
||B||^2=\frac{i}{2\pi}\int\limits_M trB\bar{B}.
$$
}
\end{theorem}
{\bf Proof.} $$ \delta_{\bar{A}}log Det \triangle_A
=\frac{\pr}{\pr s}_{|s=0}sTr(\triangle_A^{-s}
\triangle_A^{-1}\delta \triangle_A)=Tr
(\triangle_A^{-s}\bar{\pr_A}^{-1}\delta\bar{A})_{|s=0}. $$ Taking
the variation in the form
$\delta\bar{A}=\bar{\pr_A}(\epsilon)=[\bar{\pr_A},\epsilon]$ one
gets $$ Tr(\triangle_A^{-s}\bar{\pr_A}^{-1}\delta\bar{A})_{|s=0}=
Str(\triangle_A^{-s}\bar{\pr_A}^{-1}\delta\bar{A})_{|s=0}
=Str(\triangle_A^{-s}\epsilon)$$ Simple calculation shows that the
regular value of $\langle x|\triangle_A^{-s}|x\rangle$ at $s=0$ is
equal to $\frac{1}{2\pi i}(F_A + \frac{1}{2}F_{\tau_M})(x)$, where
$\tau_M$ is holomorphic tangent bundle. Thus we have $$
\delta_{\bar{A}}log Det \triangle_A = \frac{1}{2\pi
i}\int\limits_M tr F_A \epsilon + \frac{1}{2}F_{\tau_M}tr\epsilon
$$ Hence $$ \delta_A\delta_{\bar{A}}log Det \triangle_A=
-\frac{1}{2\pi i}\int\limits_M tr \bar{\pr_A}(\delta A)\epsilon =
-\frac{i}{2\pi}\int\limits_M tr\delta A\delta\bar{A}$$ Therefore
there is a holomorphic function $Det(A_0,A)$ on $\mathcal{A}$ such
that $$ Det\triangle_A=e^{-||A-A_0||^2}|Det(A_0,A)|^2 $$

 Denote $\pr= {\pr}_{A_0}$, $\bar{\pr}=\bar{\pr}_{A_0}$. Making
infinitesimal gauge transformations \bqa \delta A &=& \pr_A
(\epsilon)=\pr (\epsilon)+[A,\epsilon] \label{inf1}\\ \delta
\bar{A} &=& \bar{\pr}_A (\epsilon)=\bar{\pr}
(\epsilon)+[\bar{A},\epsilon] \label{inf2} \eqa in the formula
(\ref{quillen}) one gets $$ \delta_{\epsilon}log Det\Delta_A =
-\frac{i}{2\pi}\int tr\delta A \bar{A}-\frac{i}{2\pi}\int tr A
\delta\bar{A}+ \delta_{\epsilon}log Det (\pr
+A)+\delta_{\epsilon}log Det (\bar{\pr} +\bar{A}). $$ Since the
determinant of $\Delta_A$ is gauge invariant we have: $$
\delta_{\epsilon}log Det (\pr +A)+\delta_{\epsilon}log Det
(\bar{\pr} +\bar{A})=\frac{i}{2\pi}\int tr\pr_A (\epsilon)
\bar{A}+
   \frac{i}{2\pi}\int tr A\bar{\pr}_A (\epsilon).
$$ The simple identity $tr[A,\epsilon ]\bar{A} +tr
A[\bar{A},\epsilon]=0$ leads to $$ \delta_{\epsilon}log Det (\pr
+A)+\delta_{\epsilon}log Det (\bar{\pr}
+\bar{A})=\frac{i}{2\pi}\int tr\pr (\epsilon)  \bar{A}+
   \frac{i}{2\pi}\int tr A\bar{\pr} (\epsilon).
$$ Both left and right hand sides of the formula are decomposed
into the sum of the holomorphic and antiholomorphics parts.
Considering the antiholomorphic part we obtain the variation
formula: \bqa \label{varap} \delta_{\epsilon}log Det (\bar{\pr}
+\bar{A})=\frac{i}{2\pi}\int tr\pr (\epsilon)  \bar{A} \eqa Now
let us define the variation derivative $A(\bar{A})$ of $\log
Det(\apr +\bar{A})$ as: \bqa \label{var5}
 \delta_{\bar{A}} \log Det(\apr +\bar{A})=\frac{i}{2\pi}\int d^2z( \delta
 \bar{A}\wedge A(\bar{A})) \eqa
Expressing the equation (\ref{varap}) in terms of $A(\bar{A})$ and
$\bar{A}$ we get the condition: \bqa \label{faltap} \apr A-\pr
\bar{A} +[\bar{A},A]=0 \eqa


\end{document}